\documentclass[letterpaper,11pt]{article}
\usepackage{amsmath, amsthm, amssymb}    
\usepackage{bridges}			
\usepackage{graphicx}			
\usepackage[colorlinks=true, urlcolor=blue, citecolor=black, linkcolor=black]{hyperref}  
\usepackage{subcaption}			

\title{The {\em 4DLO} and other tubing models of $S^3$ symmetry}
\author{Chaim Goodman-Strauss\textsuperscript{1} and Eugene Sargent\textsuperscript{2}
\vspace{10pt}\\
\textsuperscript{1}National Museum of Mathematics, \\ New York, New York, U.S.A.; \\\ chaimgoodmanstrauss@gmail.com\\ 
\textsuperscript{2}Fayetteville, Arkansas, U.S.A.; \\eugene@eugenesargent.com}

\usepackage{bbm}
\def\ii{{\bf{i}}}
\def\jj{{\bf{j}}}
\def\kk{{\bf{k}}}
\def\oo{{\bf{1}}}
\def\qa{{\bf{a}}}
\def\qb{{\bf{b}}}
\def\qe{{\bf{e}}}

\def\qq{{\bf{q}}}
\def\qx{{\bf{x}}}
\def\q0{{\bf{0}}}
\def\qr{{\bf{r}}}
\def\qp{{\bf{p}}}
\def\qy{{\bf{y}}}
\def\qw{{\bf{w}}}

\begin{document}

\maketitle

\thispagestyle{empty}

\begin{abstract}
The {\em Four-dimensional Light Orchestra} or  {
\em 4DLO} was an interactive sculpture  at the National Museum of Mathematics (MoMath) from November 20, 2025 through January 2026, illustrating various sub-symmetries of the 24-cell with colored lights. This was part of a larger sequence of tubing sculptures aiming to bring to life a few lines of tables appearing in~\cite{conwayandsmith}, reprinted in~\cite{sot}, and further illuminated  in~\cite{rastanawi}. Best of all museum patrons could manipulate {\em 4DLO}'s lighting by singing and making funny noises into a microphone, and they did so with gusto. Here we describe some of the technical aspects of this sculpture and its context. 

\end{abstract}

\begin{figure}[h!tbp]
	\centering
	\includegraphics[width=2.4in]{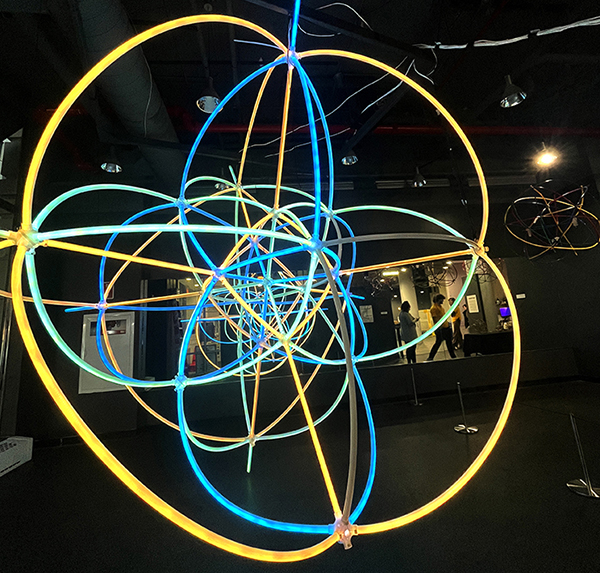} \hspace{.2in}
\includegraphics[width=3.5in]{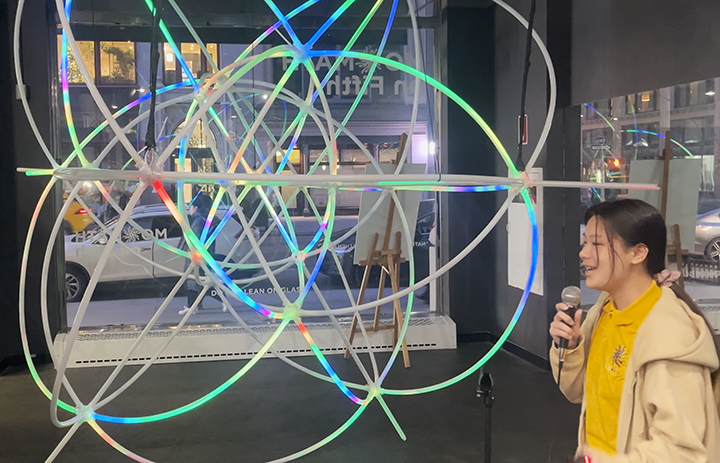} 
	\caption{(left) {\em 4DLO} showing a compound of three tesseracts, with student built work just visible in  back. (right) The colors and patterns of {\em 4DLO} responded to tone and rhythm.}
	\label{stereo}\label{fig1}
\end{figure}

\def\www{1.5in}
\begin{figure}[h!tbp]
	\centering
	(a) \ 
	\includegraphics[height=\www]{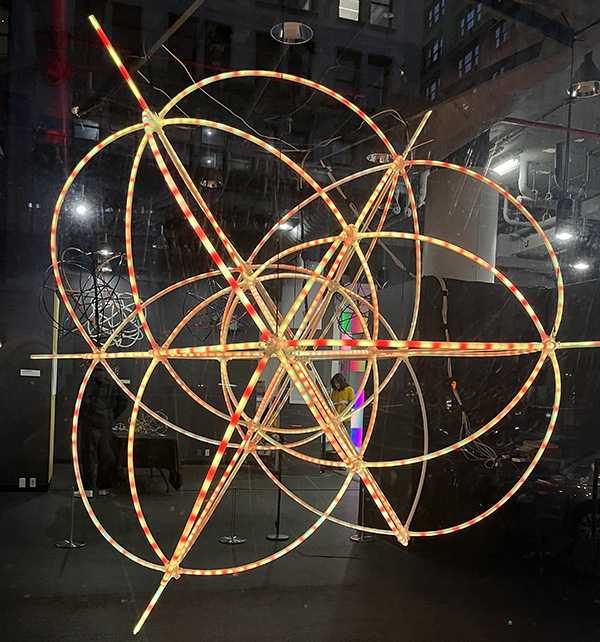} (b)\ \includegraphics[height=\www]{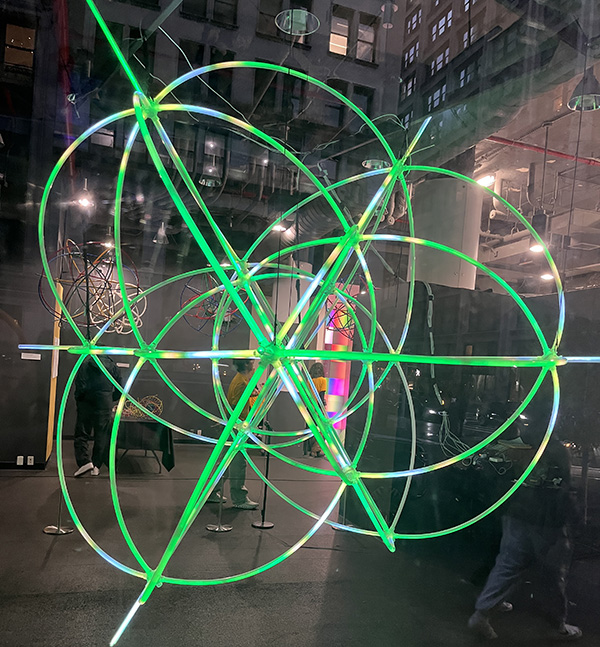}  (c)\  
\includegraphics[height=\www]{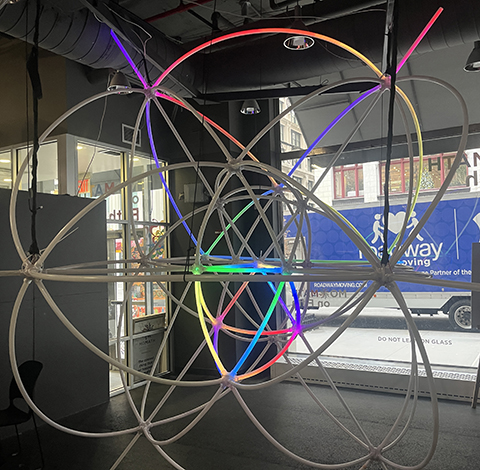}
 (d)\  \includegraphics[height=\www]{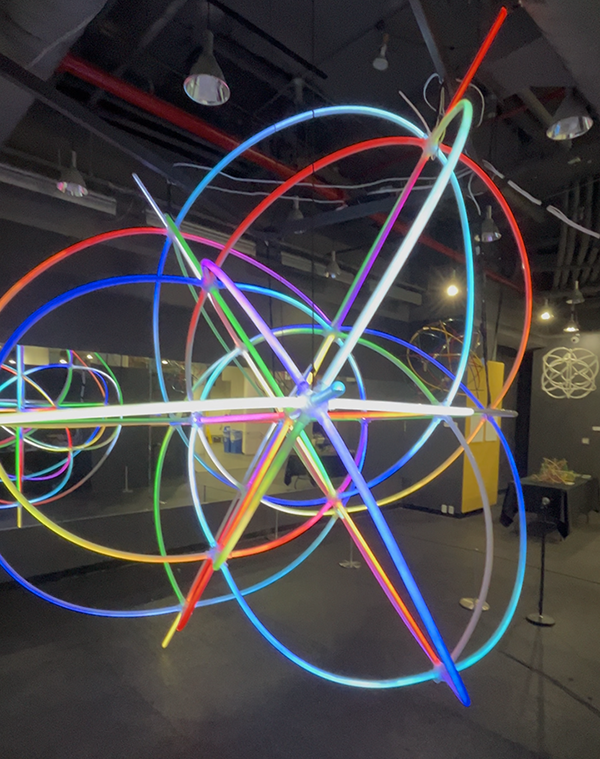}
	\caption{4DLO showing (a) the full symmetries of 24-cell; (b) half those, directing the edges; (c)  a sequence of the 
	 twenty four octahedral cells of the 24-cell, one by one (a similar sequence showed off the twenty four cubes of the three tesseracts sharing these edges); and (d) four mutually perpendicular groups of four rings of length six, indicated by color. }
	\label{hdlosequence}
\end{figure}

\section{The 24-cell and other structures in the hypersphere, $S^3$}
\vspace{1em}
The 24-cell, also known as an ``octaplex'' and by other names,  is a ``four-dimensional regular polyhedron" or {\em regular polytope}, denoted $\{3,4,3\}$ by  Ludwig Schl\"afli who  enumerated all regular polytopes in all dimensions in 1852. These are  akin to the five regular (three-dimensional) polyhedra  and to the infinitely many regular polygons in the plane.  In  each dimension $d\geq 5$  there are disappointingly only three regular polytopes, generalizations of the triangle, tetrahedron, etc., (denoted $\{3,...,3\}$ with one or more 3's); the square, cube, tesseract, etc., ($\{...,3,4\}$ with no or more 3's); and the  polytopes dual to those, a tilted square, an octahedron, the  16-cell, etc. ($\{4, 3,...\}$ with no or more 3's). In four dimensions, as it turns out, there are six. 

The list of regular polytopes was rediscovered many times.  George Hart has a helpful survey (with beautiful  3D printed models and a full bibliography) at~\cite{ghart},  H.S.M. Coxeter's {\em Regular Polytopes}~\cite{coxeter} is a canonical reference, and the Wikipedia article is extensive with detailed technical information.
  Schleimer and Segerman's 2012 Bridges paper \cite{sands} gives precise definitions of many terms we use here.

Just as a polyhedron has two-dimensional polygons for faces, and can be considered as a division of the three-dimensional sphere, a four-dimensional polytope is formed of three-dimensional cells and may be viewed as a division of the hypersphere $S^3$. In the 24-cell there  are twenty four three-dimensional octahedral cells. Each octahedron (denoted $\{3,4\}$) is composed of triangles $\{3\}$ meeting as though coming together through the sides of a square $\{4\}$. In a similar manner, at each vertex of the 24-cell, there are six octahedral cells, meeting tip to tip, as though they were coming through the square  faces of a cube $\{4,3\}$. Schl\"afli's notation encodes this data as $\{3,4,3\}$.

  This shape and further shapes derived from it have been popular, appearing over many Bridges galleries and proceedings, in various forms and to various ends:  as  collections of monkeys~\cite{Hartetal}, as a beautiful student-created aluminum frame sculpture~\cite{Peltonen}, as balloon art~\cite{Cargin}, on an engraved pendant~\cite{constant}, and as a basis for more abstract construction~\cite{Luotoniemi}. A  stainless steel 24-cell  memorial radiates in the mathematics department at Pennsylvania State University~\cite{Ocneanu}.
 
 Clearly the 24-cell is an appealing shape. Mathematically it is highly symmetrical and an outlier among the regular polytopes, with no clear analogue in other dimensions.  
 Like the tetrahedra-type polytopes  of all dimensions, the 24-cell is ``self-dual" --- at the center of each of its twenty four cells there is a vertex of another 24-cell, its dual, and {\em vice versa}.  For the same reason, it has the same number of two dimensional faces as it has edges, 96. All self-dual polytopes have palindromic Schl\"afli symbols because reversing the symbol of a regular polytope produces the symbol of its dual --- regular polygons have symbols of the form $\{n\}$ and are also self-dual. Two copies of any of these self-dual polytopes  form a natural compound.  For example, the compound of two self-dual three-dimensional tetrahedra is well-known as {\em Kepler's Stella Octangula} and two dual triangles form the familiar {\em Star of David}. So too is there a natural compound of two dual copies of the 24-cell --- a portion of this structure appears in one of Rinus Roelef's steel sculptures, landing a photo of the first author and Edmund Harriss in a newspaper~\cite{bridgesnyt}.
Among these self-dual compounds, only the 24-cell has much structure to its symmetries (that is, structure to the  subgroups of its symmetry group).

The symmetries of the hypersphere were enumerated in  stages, begun by 
 Goursat (1889), continued by Threlfall and Seifert (1931-3), Hurley (1951), Du Val (1964), and Conway and Smith (2003,~\cite{conwayandsmith}). In their recent monograph~\cite{rastanawi}, Rastanawi and Rote give a geometric re-enumeration, on which we lean heavily for evaluating  symmetry types.  (Their Section 3.3 gives a detailed survey of this history.)  
In effect, this project is to illustrate part of Conway and Smith's tables decades late for {\em The Symmetries of Things}~\cite{sot}.

\begin{figure}[h!tbp]
	\centering
\includegraphics[height=1.5in]{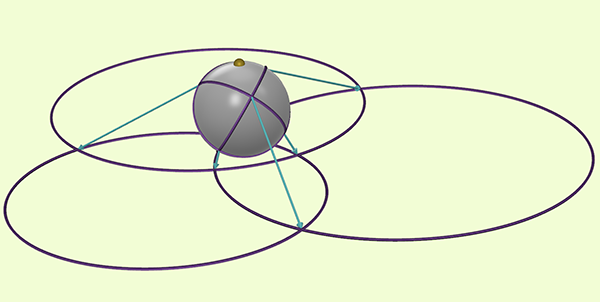}\hspace{.2in}
\includegraphics[height=1.5in]{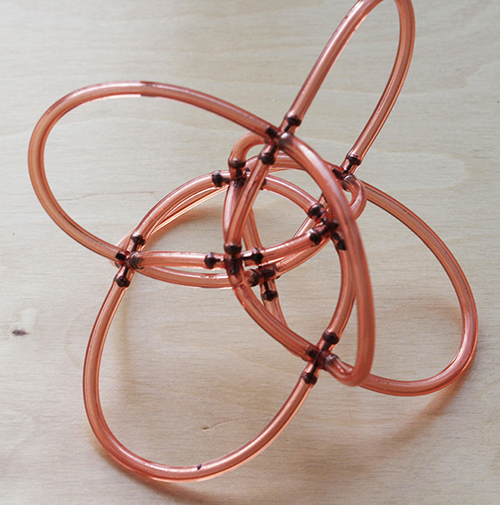}\hspace{.2in}
\includegraphics[height=1.5in]{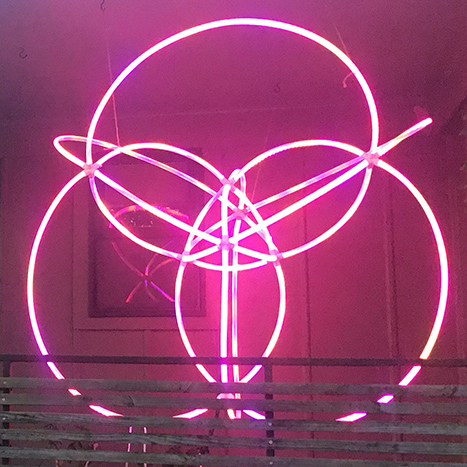}
	\caption{At left, a spherical octahedron and its stereographic projection. Both are divisions of their space into triangular cells bounded by circular arcs meeting at right angles.  Stereographic projections of  16-cells, with tetrahedral cells and faces and edges meeting at right angles: (middle) our first model made of tubing (right) a dynamic lit model in Fayetteville Ark, 2019. The polytopes in this figure are arranged ``cell-down", with a cell at the center of their projected image.}
	\label{stereo}\label{lit16}
\end{figure}

 
 \section{Rendering structures in the hypersphere}
 
The problem of visualizing the hypersphere is fundamentally three-dimensional, just as drawing pictures of features on the spherical Earth is inherently two-dimensional.  Nearly two thousand years ago Claudius Ptolemy produced the first atlas, with thousands of located sites spanning from the West African coast to the South China Sea. He projected points on Earth to a map, that is, from points on the sphere to points in a plane. Cartographers find his  conical and cylindrical projections  foundational to this day, but they've found little use for his {\em stereographic projection}, despite it being  very good for small regions away from some chosen point $N$. Mathematicians love it however, for not only does it preserve angles as measured on the sphere when they are measured on the plane,  but it also takes circles on the sphere to circles in the plane (understanding, as did Nicolaus of Cusa ca. 1450, that lines are simply circles of infinite radius). 
And this mapping  is elegantly defined: 
Geometrically, we choose some plane $P$ and a point $N$ on a sphere. Stereographic projection takes each point $p$ on the sphere to a point $\pi(p)$ in  $P$, defined as the point where the line through $p$ and $N$ meets $P$. ($\pi(N)$ is itself ``at infinity," a kind of NaN.)
 In Figure~\ref{stereo} this is shown in the standard way,  with $N$ at the top and $P$ tangent at $-N$; we say ``down" in this sense, towards this tangent. In coordinates, say $p=(x,...,z,w)$ of any dimension $d$,  it is convenient to take $N=(0,..,0,-1)$ and project to the plane with $w=0$, using 
 $\pi(p)=(x,...,z)/(1+w)$, a $(d-1)$-dimensional coordinate in $P$. 

By examining a pattern in $P$, such as the division of the plane into eight right-angled triangular regions at left in Figure~\ref{stereo}, we can understand at least something about a division of the sphere in the same way (though a spherical octahedron will have much more symmetry than any planar drawing of it).  We, as have many others, apply the same strategy one dimension up to learn something of the possible symmetries in the hypersphere, like those of the  ``16-cells" that are shown middle and right in Figure~\ref{stereo}.

\def\re{{\text{Re}}}
\def\im{{\text{Im}}}

Tubing is manufactured in circular coils, and therefore pieces of it relax into circular arcs --- this insight, that tubing is great for stereographic renderings of circles in the hypersphere, was the aha moment that spurred on this entire multifaceted project.  With no delay the first author sought what might be at hand, finding some vinyl tubing at an automotive store and six-way connectors as playing jacks at a party store. The first of these models, a 16-cell, is at middle in Figure~\ref{stereo}, and the second is  at left  in  Figure~\ref{16cellsequence}.

 
 \begin{figure}[h!tbp]
	\centering
	
\includegraphics[height=2in]{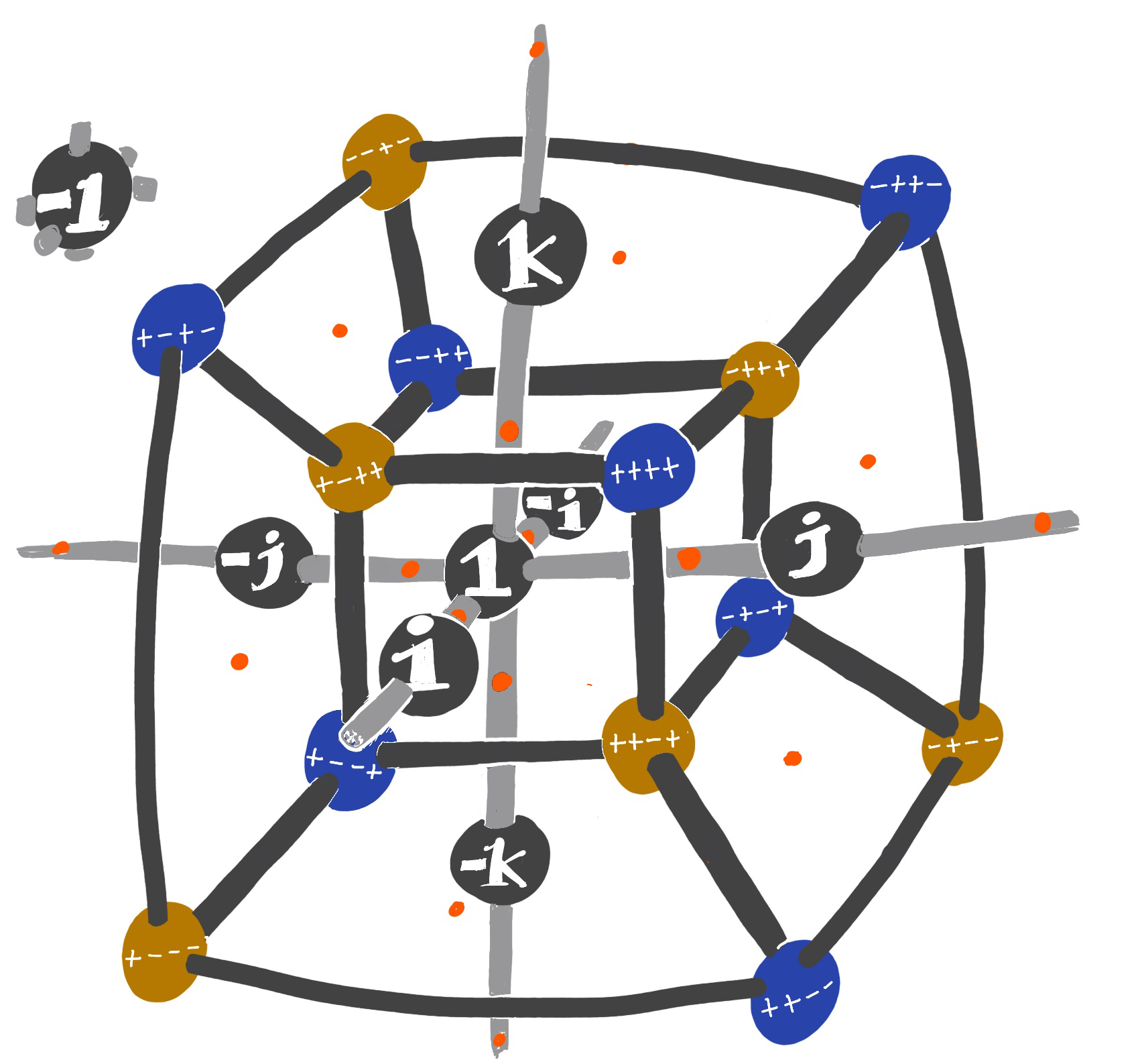}
\hspace{.2in}
\includegraphics[height=2in]{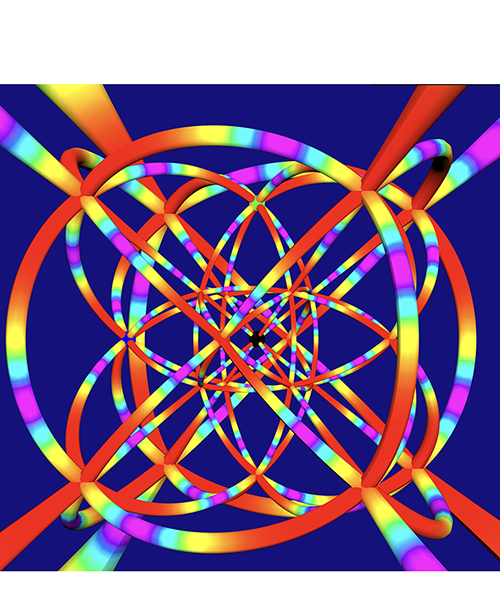}

\caption{(Left) With some deletions and distortions, stereographically projected sets $V_8$ (black), $V_{16}^+$ (blue), and $V_{16}^{-}$ (brown); with edges of a cell-down tesseract (black); the red points on the faces of this tesseract are the points in $V_{24}^-$. It can help to sketch your own version of this diagram. (Right) a screenshot of our interactive javascript model is at {\em chaimgoodmanstrauss.com/4DLO}.}\label{coords}
\end{figure}

\begin{figure}[h!tbp]
	\centering
\includegraphics[height=1.8in]{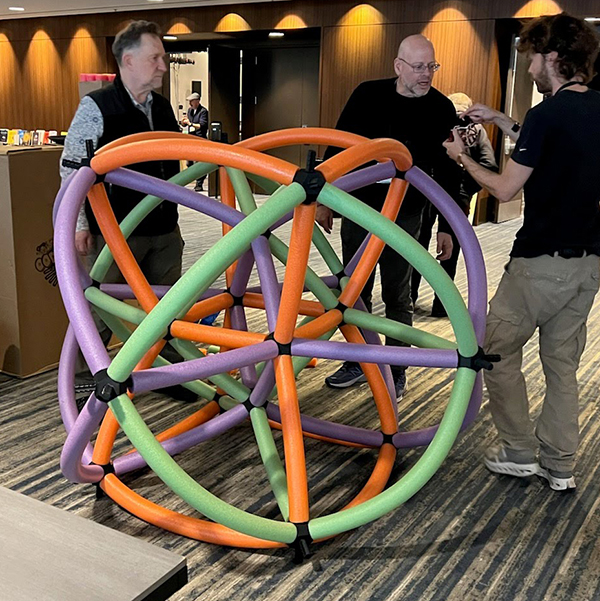}
\hspace{.3in}
\includegraphics[height=1.8in]{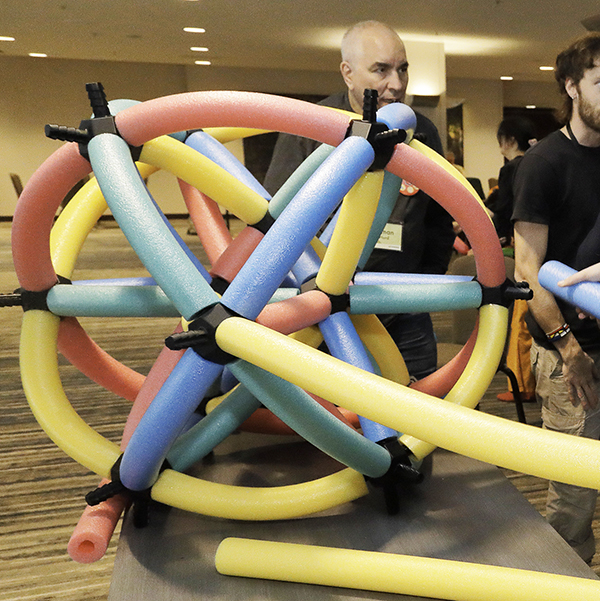}
\caption{Cell-down  24-cells, shown as (left)  a compound of three tesseracts and (right)  $T^*$ generated by the three-fold rotational symmetries of a tetrahedron. As with many of these models, not all edges of the mathematical structure are shown. The models in this figure are missing edges outside a unit ball, the ``top" half of the 24-cell. Scott Vorthmann led a group build of these  models at the sixteenth Gathering for Gardner, G4G16 (Feb. 2026).}\label{celldown24s}\end{figure}

\begin{figure}[h!tbp]
	\centering
	\includegraphics[width=\textwidth]{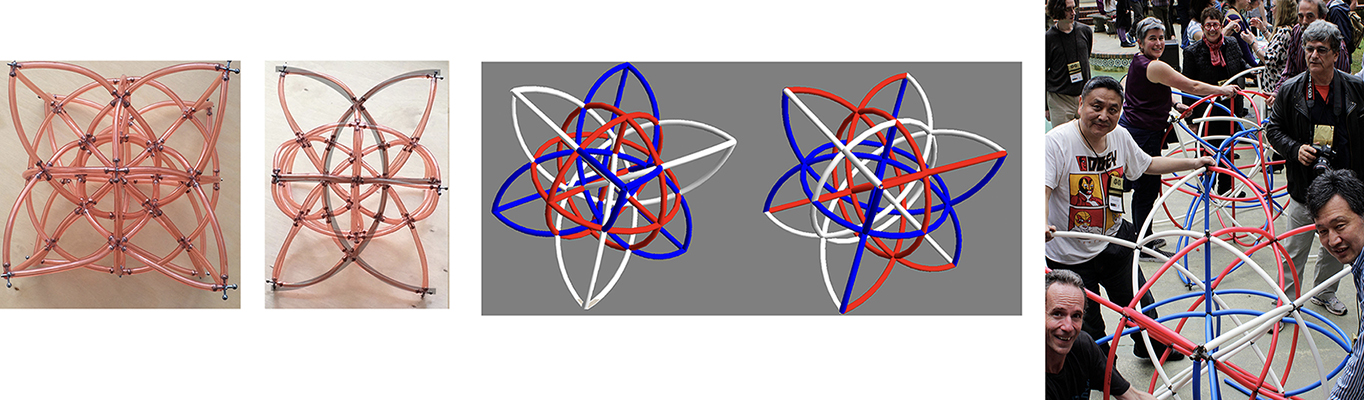}

	\caption{(At left) A first model of a compound of two compounds of three 16-cells (
	``two$\cdot$ three 16-cells'')  made of  playing jacks and  automotive tubing; the marked half-circles lie in a  $\sqrt 2$-to-1 rectangle (see Fig.~\ref{plans}).
	 (Middle) Six 16-cells in two sets of three: one vertex-down in red, with vertices $V_8$; $V_{16}^+$ and $V_{16}^-$ in white and blue; the other with vertices in $V_{24}'$. 
	(Right) {\em A Compound of Two Compounds of Three Sixteen Cells}, G4G12 (2016).}
	\label{16cellsequence}
\end{figure}

\begin{figure}[h!tbp]
	\centering
\includegraphics[height=1.4in]{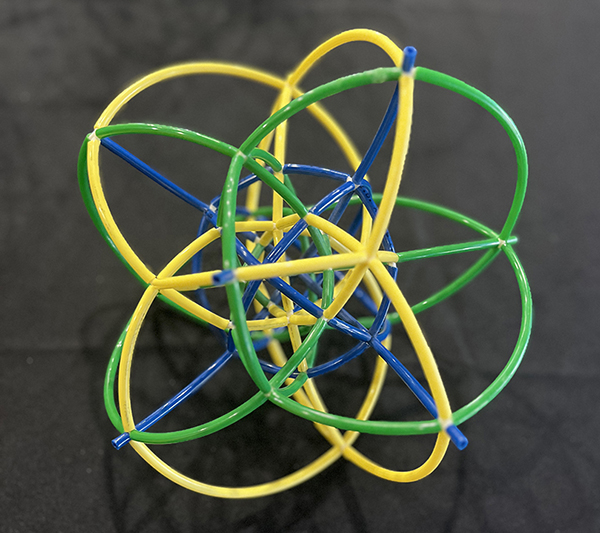}
\hspace{.3in}
\includegraphics[height=1.4in]{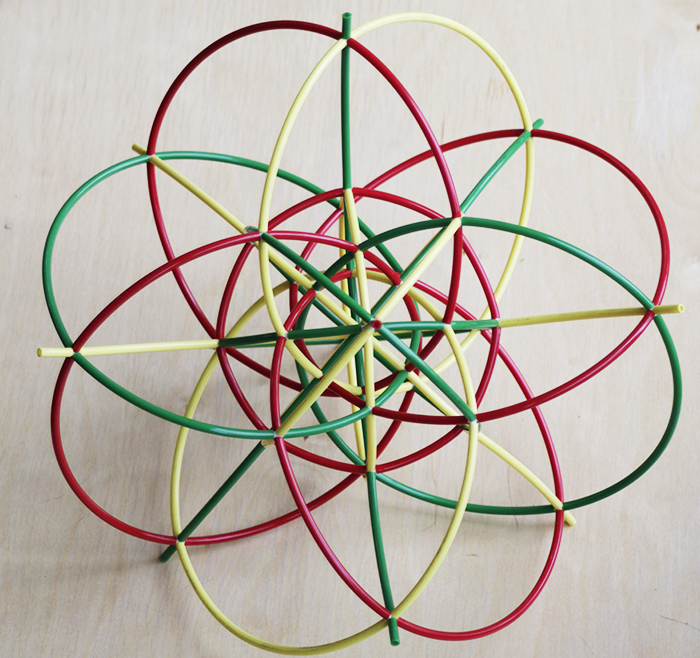}
\hspace{.3in}
 \includegraphics[height=1.4in]{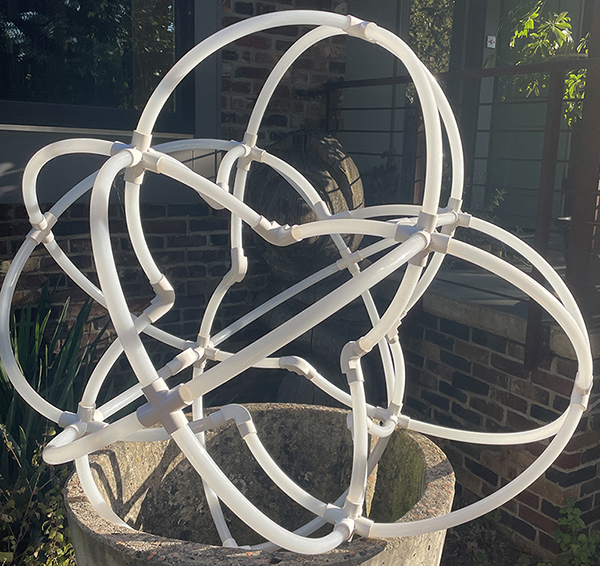}
\caption{(Left) an obliquely arranged  compound of three 16-cells.
(Middle) if a direction were added to these edges, this model would show symmetry $O^*$. (Right) a pleasing piece of yard art, a portion of three 16-cells made from four-way and two-way connectors.}
\label{three16cellstimes2}
\end{figure}

 We describe coordinates as quaternions, indicating $(a,b,c,d)$ as $\qq=a\ii+b\jj+c\kk+d\oo$. On these we allow the usual rules of arithmetic, with multiplicative identity $\oo$, plus the rules that $\ii\jj=\kk$ and, for any {\em imaginary} (that is,  $\re(\qq):=d$ is zero) and  {\em unit} ($a^2+b^2+c^2+d^2=1$) quaternion $\qq$, we have $\qq\qq=-\oo$.  (See~\cite{sands} for further details.) The unit quaternions are exactly the points of the hypersphere $S^3$, and therefore coordinatize it. For illustrations, stereographically projecting from $-\oo$ into our space places $\oo$ at the origin and $-\oo$ ``at infinity", the unit imaginary quaternions on the unit sphere and $\pm\ii, \pm\jj$ and $\pm\kk$ upon it as coordinate axes. Our sketch is at left in Figure~\ref{coords} and the active reader is encouraged to draw their own. The familiar ``cell-down" tesseract $\{4,3,3\}$ sits cell-downwards in $R^4$. It has sixteen vertices, we'll call $V_{16}$, at points of the form $\frac{1}{2}(\pm \ii \pm \jj \pm \kk \pm \oo)$. In stereographic projection these lie at the corners of two nested cubes to either side ($\pm\oo$) of the unit sphere (and a cell is at the center of the projection). We'll further sort these points into the ``even" ones $V_{16}^+$ with an even number of $-$'s, and its ``odd" complement $V_{16}^-$. The edges of a tesseract connect vertices that differ by exactly one coordinate, and therefore each edge has one end in 
  $V_{16}^+$ and the other in $V_{16}^-$.    Three tesseracts can be arranged as at left in Figure~\ref{fig1};  the yellow one, though tilted to the floor, is cell-down (a cube is at its center) and has  vertices in $V_{16}$. 
It is dual to a ``vertex-down" 16-cell with eight vertices, $V_8=\{\pm\oo, \pm\ii,\pm\jj,\pm\kk\}$ (One appears in red at the center of Figure~\ref{16cellsequence}; you can quickly sketch a vertex-down 16-cell by connecting up the points in $V_8$ on your drawing.)  ``Cell-down" 16-cells with vertices in $V_{16}^+$ are shown in Figure~\ref{stereo}. In fact  $V_{16}\cup V_8$, which we denote $V_{24}$, form the vertices a 24-cell, which we'll call  ``vertex-down";  its dual ``cell-down"  24-cell  has vertices $V_{24}'$ of the form $\frac1{\sqrt 2}(\pm\qx\pm\qy)$ where $\qx, \qy$ are two of    
    $\oo, \ii,\jj$, and $\kk$. 
    The {\em 4DLO} of Figures~\ref{fig1} and~\ref{hdlosequence} show a vertex-down 24-cell (and thus with a vertex at the center of the projection), but with coordinates tilted so that  gravity pulls in the $(-1,0,-1)$ direction.  The 24-cells of  Figure~\ref{celldown24s},  made at the sixteenth Gathering For Gardner  in February, 2026,  are cell-down. 
    
The vertex-down 24-cell shares its vertices $V_{24}$ with a compound of three 16-cells, with vertices $V_8, V_{16}^+$ and $V_{16}^-$, and thus there is another compound of three 16-cells with corresponding vertices in $V_{24}'$.  We show these in Figures~\ref{16cellsequence} and~\ref{three16cellstimes2}.
   In each compound of three 16-cells, the vertices of one 24-cell are those of the 16-cells, but the vertices of the dual 24-cell are  merely points where the edges of each 16-cell appear to cross, just as  edges of dual tetrahedra cross in a stella octangula.  The two compounds of three 16-cells use the same edges but swap the role of the vertices! (The model at left in Figure~\ref{16cellsequence} is both  together overlapping.) Together these two compounds form 
         {\em A Compound of Two Compounds of Three 16-cells}, assembled  in two physical pieces in Atlanta in 2016 by 
 participants of the twelfth biennial Gathering for Gardner,  one of a sequence of sculptures we were fortunate to have created under the inspiring, supportive and generous leadership of George Hart, familiar to many of us in the math-art community. 
 
\begin{figure}[h!tbp]
	\centering
\includegraphics[height=1.4in]{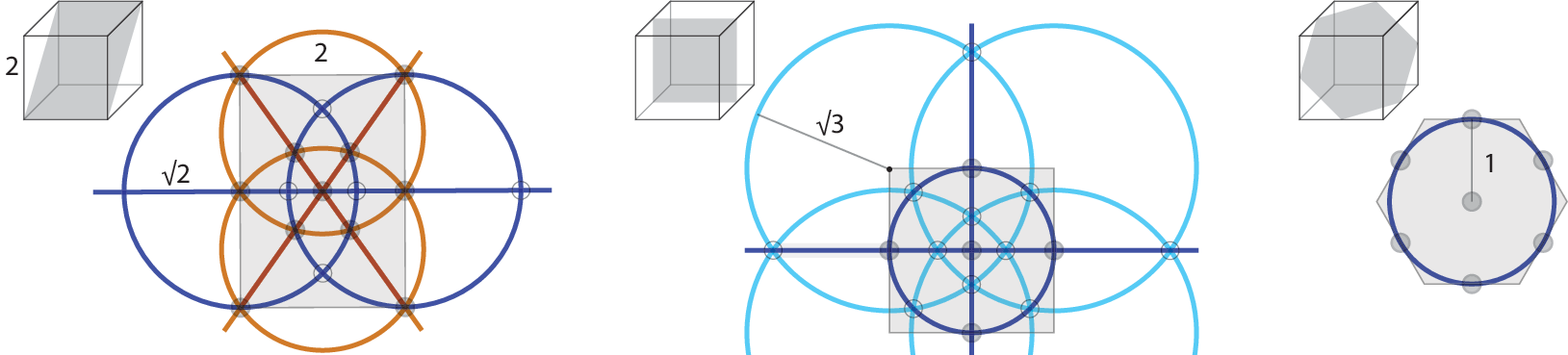}
	\caption{Plans for a vertex-down projection of the 24-cell (brown);  a cell-down 24-cell (light blue); and a compound of two compounds of three 16-cells (black); The points of $V_{24}$ are shaded and those of $V_{24}'$ circled. These plans are on planes meeting the unit cube (of width 2) as shown.}
	\label{plans}
\end{figure}

 There is a very nice compound of three tesseracts, dual to the compound of three 16-cells, sharing vertices $V_{24}$, which we show at left in Figures~\ref{fig1} and~\ref{celldown24s}. (We don't show the natural compound of six tesseracts.) You can find plans for making your own models in Figure~\ref{plans}.

\section{Quaternionic groups as the symmetries of $S^3$}

Nicely, quaternionic multiplication can describe the symmetries of  objects in $S^3$. Each given unit quaternion $\qq\in S^3$ also {\em acts} upon on $S^3$ as a whole, sliding everything together along a twisting family of great-circle rails (a ``Hopf fibration''), just by multiplying every point by $\qq$. However multiplying on the left ($\qx\mapsto\qq\qx$) slides along rails in right-handed  [sic] coils and multiplying on the right ($\qx\mapsto\qx\qq$) slides along left-handed coils (as measured against a right-handed ordering of $\ii$, $\jj$, and $\kk$). We can think of these coiling around two special ``straight" rails, 
 the one  
passing through $\oo$ and $\qq$, and its opposite, the one furthest one away and perpendicular to it. (These are special only in our coordinates: geometrically every rail has an  opposite and any opposite pair of rails can be rotated to any other.) As $S^3$ slides  one full circuit {\em along} a special rail, it  makes a full revolution {\em around} it (the same thing as making a full circuit {\em along} its opposite rail).

Or we may specify a fibration of $S^3$ into great circle rails by giving some unit $\qq$ on the rail through $\oo$. We can choose that  $\qq$ to be  imaginary, and therefore on the unit great sphere midway between  $ \oo$ and $-\oo$. Multiplying by $\qr := (\cos\theta) \oo+(\sin\theta) \qq$ moves the hypersphere $\frac{\theta}{2\pi}$th the way along this rail, but twisting with a handedness: Right multiplication $\qx\mapsto\qx\qr$ rotates $S^3$ by $\theta$ clockwise, as seen from $\qr$ looking to $\oo$; and left-multiplication  $\qx\mapsto\qr\qx$ rotates $S^3$ by $\theta$ counter-clockwise. Therefore taking $\qx\mapsto{\qr}^{-1}\qx\qr$ slides $S^3$ along, then back, fixing $\oo$ but rotating $S^3$ about the $\oo$-to-$\qr$ axis clockwise by $2\theta=2\arccos(\re(\qr))$. We can thusly encode three-dimensional rotations by unit quaternions $\qr$. The correspondence is two-to-one, though, because ${\qr}^{-1}\qx\qr={(-\qr)}^{-1}\qx(-\qr)$ and so $\qr$ and $-\qr$ give the same rotation of $\oo$.
Each 3D rotational symmetry group $G$ lifts to a symmetrical set of points $G^*$ in $S^3$ (in binary, opposite pairs).  Multiplying on the right (say) by any of these points $\qq$ rotates $S^3$ taking $\oo$ to $\qq$ and preserving $G^*$ as a whole --- $G^*$ is a fixed-point free symmetry group of $G^*$ itself. V. Hart {\em et al.} cleverly  show  these ``binary'' rotation groups of the platonic solids using monkeys~\cite{Hartetal}.

 The 24-cell is related to the symmetries of the tetrahedron in this way (and a dual pair of 24-cells is related to the symmetries of a cube). To see this orient, a tetrahedron so that its vertices are $(1,1,1)$ and permutations of  $(-1,-1,1)$ (that is, alternating vertices of the unit cube),  and let $T$ be its rotational symmetry group. 
 The four pairs of antipodal points in $V_8$ correspond to the identity and the three 2-fold rotations, and the eight pairs of antipodal points in $V_{16}$ correspond to the eight 3-fold rotations. In fact, edges correspond to rotations:  If $\qp,\qr\in V_{24}$ are the ends of an edge  then $\qe := (\qp\qr^{-1})$ corresponds to a three-fold rotation (that is $\qe ^3=-\oo$), and {\em vice versa}. Moreover rotations that have the same axis correspond to edges lying on the same great circle paths. 
 Figure~\ref{hdlosequence}(d) shows this on {\em 4DLO}.  At right in Figure~\ref{celldown24s} we color cycles of 24-cell edges by which three-fold rotation axis they correspond to, showing symmetry $T^*$. 

We can map the 24-cell to its dual, taking $V_{24}$ to $V_{24}'$ by multiplying (right or left) by $\frac{1}{\sqrt{2}}(\oo+\ii)$, or by any other point in $V_{24}'$. The twelve pairs of points in $V_{24}'$  each  correspond to a quarter turn  about a coordinate axis, exactly the remaining rotational symmetries $O$ of a cube. The binary octahedral group $O^*=V_{24}\cup V_{24}'$. Connecting these vertices by edges corresponding to four-fold rotations we have the model at  left in Figure~\ref{16cellsequence}, or if we color these by rotation axis, we have the one at right in Figure~\ref{moremodels}. 
Directing these edges gives symmetry $O^*$, but these models have more symmetry,  because we can pivot around points: By multiplying on both left and right, we can fix points, and therefore have symmetry about them. 
 \def\www{1.4in}
\begin{figure}[h!tbp]
	\centering
 \includegraphics[height=\www]{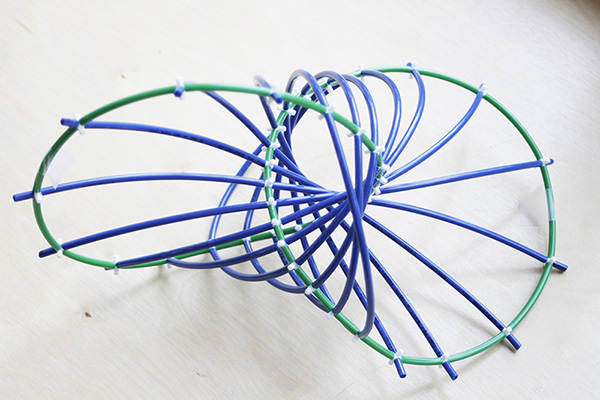}\hspace{.3in}
\includegraphics[height=\www]{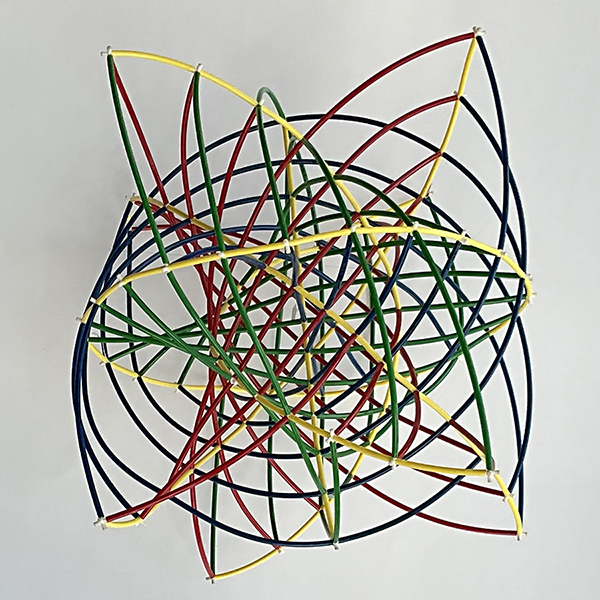}\hspace{.3in}
\includegraphics[height=\www]{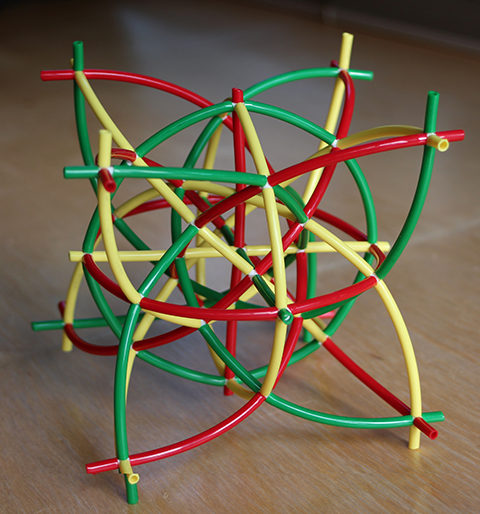}
	\caption{Models of symmetry types  (left) $\pm[C_2\times C_{11}]$,   (middle) $\pm[T\times C_{6}]$, six copies of a 16-cell slid along yellow rails, and (right) $O^*$ with edges corresponding to four-fold rotations of the cube.}
	\label{moremodels}
\end{figure}
 By writing $\pm[A\times B]$ Conway and Smith~\cite{conwayandsmith} mean the group of  actions of the form  $\qx\mapsto \qa\qx \qb$ where the quaternions $\qa$ and $\qb$ correspond to rotations that lie in symmetry groups $A$ and $B$. 
 Figure~\ref{moremodels} shows models of symmetry $\pm[C_2\times C_{11}]$ and $\pm[T\times C_{6}]$, where $C_n$ is the cyclic group of order $n$. 
 That notation soon has many amendments ---  we rely  on the detailed geometric discussion in ~\cite{rastanawi}, especially their Sections 8.4 and 8.6 on the tesseract and 24-cell.   
 We'll use Conway and Smith's notation, referring the reader to the concordances and extensions in \cite{conwayandsmith, sot, rastanawi}.

 As there is cubic symmetry around each point in $V_{24}\cup V_{24}'$, this compound of two 24-cells  has rotational symmetry $\pm[O\times O]$.   Applying any $\qx\mapsto \qr\qx\qp = (\qr\qp)(\qp^{-1}\qx\qp)$ rotates $\oo$ about the $\qp$ axis (a rotation fixing $\oo$ that is a symmetry in $O$) and then by left multiplication moves everything taking $\oo$ to $\qr\qp$, another point of $V_{24}\cup V_{24}'$, also in $O^*$.
The model at  left in Figure~\ref{16cellsequence} also has vertices $V_{24}\cup V_{24}'$, and has this same symmetry, but it is not two dual 24-cells.
Allowing mirror reflections by taking symmetries of the form $\qx\mapsto \qp\overline{\qx}\qr$ where
 $\overline{\qx}=2\re(\qx)-\qx$, we extend these to the full symmetry of this dual pair,  named 
 $\pm[O\times O]\cdot 2$. Henceforth we'll name only the rotational symmetries.
In  Figure~\ref{hdlosequence}(a) 4DLO shows the full symmetry of a single 24-cell   
denoted  $\pm\frac12 [O\times O]$, with half  the symmetry of a dual pair.  We've seen that the vertices $V_{24}$ of a vertex-down 24-cell sort into the sets $V_{8}$, $V_{16}^+$ and $V_{16}^-$.    Any point corresponding to 3-fold rotation, so it happens, is a  symmetry that permutes these sets.  For example, multiplying by $\qw:= \frac12(\ii+\jj+\kk+\oo)$ moves $\oo\in V_{8}$ to 
$\qw\in V_{16}^+$, to $\qw^2\in V_{16}^-$, to $-\oo\in V_8$, to $-\qw\in V_{16}^+$, etc. {We can put a direction on the edges of a 24-cell in this way, with type (we think)  $\pm[T\times T]$, shown at right in Figure~\ref{coords} and in Figure~\ref{hdlosequence}(b).}

Because $V_{8}$ is the vertex set of a 16-cell, its images $V_{16}^+$ and $V_{16}^-$ are too,  explaining the symmetry in the compound of three 16-cells with vertices $V_{24}$, and in the  self-dual compound of two of those compounds,  shown in Figures~\ref{16cellsequence} and~\ref{three16cellstimes2}.  Equivalently, the 16-cell, as well as its dual, the tesseract, has a third of the symmetry of a 24-cell, which Conway and Smith list as $\pm\frac16[O\times O]$.   Dually there is compound of three tesseracts sharing vertices with a (dual) 24-cell, shown cell-down in 
 Figure~\ref{celldown24s}, or  vertex down in {\em 4DLO} at left in Figure~\ref{fig1}, and a compound of two compounds of  three tesseracts, with the full symmetry $\pm[O\times O]$.

%

\section{{The Four-Dimensional Light Orchestra}, a.k.a. {\em 4DLO}}

Inspired by Leo Villareal's {\em Buckyball} and {\em Hive: Bleeker Street}, we turned to moving  lights to illustrate these symmetries. A  first project is shown in the center in Figure~\ref{lit16}. The {\em 4DLO} was much more ambitious, consisting of about 300 feet of  1 inch outer diameter polyethylene tubing (with its beautiful capacity for radiating soft light), and 14,000 WS2811b LED lights diffused by strips of packing foam. 
These are controlled very simply -- each light on each strand receives a string of bytes, taking three for its R, G, and B values and passing along the rest.  The 96 tubes of the model were organized into four quadrants, each with seven separate strands, powered at each node they passed through. A 
Teensy 4.1 microcontroller has more than enough processing power and pins for this task, especially when decked out with a full suite of accessories from its producer PJRC Electronic Products: four Octos, one for each quarter for the data lines, an Audio Shield for sound, and  extra on-board memory. Each individual LED was associated with a particular location in 4D space, measured along the arc of an edge of the 24-cell;  at each moment of the program, an edge would be associated with a coloring function, and these would be assigned by using whatever symmetry we wished to highlight. The javascript program at {\em https://chaimgoodmanstrauss.com/4DLO/} allowed us to prototype  the geometry and symmetry --- there you can also rotate the 24-cell into different positions. We used a bespoke scripting language for  the sequencing, which we developed making extensive use of  {\em Claude.ai} (along the way, learning much about the potential and limitations of today's LLMs).  All of our code is at the Github repository {\em github.com/chaimgoodmanstrauss/4DLO.git}.
Needless to say, we made many mistakes, and  learned a lot, much too much to fit here. (We printed four  complete sets of nodes before getting it right.) Give us a call if you would like to know more. 


%



\section*{Acknowledgements} We thank the referees for improving this paper  with  many helpful comments.
Figure~\ref{celldown24s}  courtesy (left) National Museum of Mathematics (MoMath), (right)  Gathering For Gardner Foundation (G4G), Stephen Inglima.
This work has been supported by G4G; MoMath; Institut Henri Poincar\'e (UAR 839 CNRS-Sorbonne Universit\'e, Labex Carmin ANR-10-LABX-59-01');  and Universidad Nacional Autónoma de México. We  are especially grateful to unnamed donors, to Cindy Lawrence, and to the MoMath staff for helping bringing {\em 4DLO} to life. 
The first author expresses deep appreciation for the second author's ability to make these things real. 
    
{\setlength{\baselineskip}{13pt} 
\raggedright				

} 
   
\end{document}